\documentclass{gtart}

\def\ifplaintex{\expandafter\ifx\csname documentclass\endcsname\relax}

\ifplaintex 
\else
\fi

\expandafter\ifx\csname beginpicture\endcsname\relax
\expandafter\ifx\csname documentclass\endcsname\relax
\input pictex \else
\input prepictex \input pictex \input postpictex \fi\fi

\def\gt{{\mathsurround=0pt\it $\cal G\mskip-2mu$eometry \&\ 
$\cal T\!\!$opology}}        

\def\gtp{{\mathsurround=0pt\it $\cal G\mskip-2mu$eometry \&\ 
$\cal T\!\!$opology $\cal P\!$ublications}}  

\def\lognumber#1{\def\thelognumber{#1}}
\def\volumenumber#1{\def\thevolumenumber{#1}}
\def\papernumber#1{\def\thepapernumber{#1}}
\def\volumeyear#1{\def\thevolumeyear{#1}}

\def\pagenumbers#1#2{\def\startpage{#1}\def\finishpage{#2}}
\def\published#1{\def\publishdate{#1}}
\def\proposed#1{\def\theproposer{#1}}
\def\seconded#1{\def\theseconders{#1}}
\def\received#1{\def\receiveddate{#1}}
\def\revised#1{\def\reviseddate{#1}}
\def\accepted#1{\def\accepteddate{#1}}
\def\asciititle#1{\def\theasciititle{#1}}

\def\asciikeywords#1{\def\theasciikeywords{#1}}

\def\shorttitle#1{\def\theshorttitle{#1}}

\let\\\par\let\thelognumber\relax
\let\thevolumenumber\relax\let\thepapernumber\relax
\let\thevolumeyear\relax\let\thesamplenumber\relax\let\startpage\relax
\let\finishpage\relax\let\publishdate\relax\let\receiveddate\relax
\let\reviseddate\relax\let\accepteddate\relax\let\theasciititle\relax
\let\theasciiauthors\relax
\let\theasciiabstract\relax\let\theasciikeywords\relax
\let\theasciiemail\relax\let\theshortauthors\relax\let\theshorttitle\relax

\long\def\maketitlep{   

\count0=\startpage

\gt\hfill      
\beginpicture
\setcoordinatesystem units <0.33truein, 0.33truein> point at 2.2 0.9
\setplotsymbol ({$\cal G$})
\plotsymbolspacing=9truept
\circulararc 315 degrees from 0 1 center at 0 0
\setplotsymbol ({$\cal T$})
\circulararc 315 degrees from 1 -1 center at 1 0
\endpicture
\break
{\small\ifx\thesamplenumber\relax 
Volume \else Sample
\fi\thevolumenumber\ (\thevolumeyear)
\startpage--\finishpage\nl
Published: \publishdate}
\vglue 0.5truein plus 0.4fil minus 0.1truein

{\parskip=0pt\leftskip 0pt plus 1fil\def\\{\par\smallskip}{\ifplaintex\large
\else\Large\fi\bf\thetitle}\par\medskip}   

\vglue 0pt plus 0.1fil 

{\parskip=0pt\leftskip 0pt plus 1fil\def\\{\par}{\sc\theauthors}
\par\medskip}

\vglue 0pt plus 0.1fil 

{\small\parskip=0pt\let\newline\\
{\leftskip 0pt plus 1fil\def\\{\par}{\sl\theaddress}\par}
\expandafter\ifx\theemail\relax    
\relax\else\vglue 5pt plus 0.02fil minus 2pt\def\\{\stdspace{\rm 
and}\stdspace} 
\cl{Email:\stdspace\tt\theemail}\fi
\ifx\theurl\relax                  
\relax\else\vglue 5pt plus 0.02fil minus 2pt\def\\{\stdspace{\rm 
and}\stdspace}
\cl{URL:\stdspace\tt\theurl}\fi\par}

\vglue 7pt plus 0.3fil minus 3pt

{\bf Abstract}
\vglue 5pt plus 0.1fil minus 2pt

\theabstract

\vglue 7pt plus 0.3fil minus 3pt

{\bf AMS Classification numbers}\quad Primary:\quad \theprimaryclass

Secondary:\quad \thesecondaryclass

\vglue 5pt plus 0.3fil minus 2pt

{\bf Keywords:}\quad \thekeywords

\vglue 10pt plus 0.5fil minus 5pt

{\small  Proposed: \theproposer\hfill Received: \receiveddate\nl
Seconded: \theseconders\hfill 
\ifx\reviseddate\relax                         
Accepted: \accepteddate                        
\else
Revised: \reviseddate                          
\fi}
\eject
}       

\let\maketitlepage\maketitlep
\let\maketitle\maketitlepage

\font\phead=cmsl9 scaled 950
\font=cmsl9 scaled 1050
\font\pnum=cmbx10 scaled 913
\font\lnum=cmbx10 
\font\pfoot=cmsl9 scaled 950
\font=cmsl9 scaled 1050
\ifplaintex
\headline{\vbox to 0pt{\vskip -4.5mm\line{\small\phead\ifnum
\count0=\startpage ISSN 1364-0380 (on line)
1465-3060 (printed) \hfill {\pnum\folio}\else\ifodd\count0\def\\{ }%
\ifx\theshorttitle\relax\thetitle\else\theshorttitle\fi\hfill{\pnum\folio}
\else\def\\{ and }{\pnum\folio}\hfill\ifx\theshortauthors\relax\theauthors
\else\theshortauthors\fi\fi\fi}\vss}}
\footline{\vbox to 0pt{\vglue 0mm\line{\small\pfoot\ifnum\count0=\startpage
\copyright\ \gtp\hfill\else
\gt, Volume \thevolumenumber\ (\thevolumeyear)\hfill\fi}\vss
}}
\else
\makeatletter
\def\@oddhead{{\small\lhead\ifnum\count0=\startpage ISSN 1364-0380 (on line)
1465-3060 (printed) \hfill {\lnum\number\count0}\else\ifodd\count0
\def\\{ }\ifx\theshorttitle\relax \thetitle \else\theshorttitle\fi\hfill
{\lnum\number\count0}\else\def\\{ and }{\lnum\number\count0}
\hfill\ifx\theshortauthors\relax 
\theauthors\else\theshortauthors\fi\fi\fi}}\def\@evenhead{\@oddhead}
\def\@oddfoot{\small\lfoot\ifnum\count0=\startpage\copyright\ \gtp\hfill\else
\gt, Volume \thevolumenumber\ (\thevolumeyear)\hfill\fi}
\def\@evenfoot{\@oddfoot}
\makeatother
\fi

\newwrite\gtoutfile
\long\gdef\makeheadfile{  
{\def\\{, }\def\s{ }
\immediate\openout\gtoutfile head.xxx
\immediate\write\gtoutfile{To: math@arxiv.org}
\immediate\write\gtoutfile{Subject: put or rep NNNNN:pppp}
\immediate\write\gtoutfile{--text follows this line--}
\immediate\write\gtoutfile{Proxy-for: \ifx\theasciiauthors\relax
\theauthors\else\theasciiauthors\fi\s<\ifx\theasciiemail\relax\theemail\else\theasciiemail\fi>}
\immediate\write\gtoutfile{\noexpand\\}
\immediate\write\gtoutfile{Authors: \ifx\theasciiauthors\relax
\theauthors\else\theasciiauthors\fi}
\immediate\write\gtoutfile{Title: \ifx\theasciititle\relax
\thetitle\else\theasciititle\fi}
\immediate\write\gtoutfile{Subj-class: GT or SG or MG etc}
\immediate\write\gtoutfile{MSC-class: \theprimaryclass\ifx\thesecondaryclass\relax\else, \thesecondaryclass\fi}
\immediate\write\gtoutfile{Journal-ref: Geom. Topol. \thevolumenumber
(\thevolumeyear) \startpage-\finishpage}
\immediate\write\gtoutfile{Comments: Published by Geometry and Topology at}
\immediate\write\gtoutfile{\s\s http://www.maths.warwick.ac.uk/gt/GTVol\thevolumenumber/paper\thepapernumber.abs.html}
\immediate\write\gtoutfile{\noexpand\\}
\immediate\write\gtoutfile{}
\ifx\theasciiabstract\relax
\immediate\write\gtoutfile{\theabstract}\else
\immediate\write\gtoutfile{\theasciiabstract}\fi
\immediate\write\gtoutfile{}
\immediate\write\gtoutfile{\noexpand\\}
\immediate\write\gtoutfile{}
\immediate\closeout\gtoutfile}}  

\def\maketitlepage{\maketitlep\makeheadfile}
\let\maketitle\maketitlepage

\def\ifplaintex{\expandafter\ifx\csname documentclass\endcsname\relax}

\ifplaintex 
\else
\fi

\expandafter\ifx\csname beginpicture\endcsname\relax
\expandafter\ifx\csname documentclass\endcsname\relax
\input pictex \else
\input prepictex \input pictex \input postpictex \fi\fi

\def\gt{{\mathsurround=0pt\it $\cal G\mskip-2mu$eometry \&\ 
$\cal T\!\!$opology}}        

\def\gtp{{\mathsurround=0pt\it $\cal G\mskip-2mu$eometry \&\ 
$\cal T\!\!$opology $\cal P\!$ublications}}  

\def\lognumber#1{\def\thelognumber{#1}}
\def\volumenumber#1{\def\thevolumenumber{#1}}
\def\papernumber#1{\def\thepapernumber{#1}}
\def\volumeyear#1{\def\thevolumeyear{#1}}

\def\pagenumbers#1#2{\def\startpage{#1}\def\finishpage{#2}}
\def\published#1{\def\publishdate{#1}}
\def\proposed#1{\def\theproposer{#1}}
\def\seconded#1{\def\theseconders{#1}}
\def\received#1{\def\receiveddate{#1}}
\def\revised#1{\def\reviseddate{#1}}
\def\accepted#1{\def\accepteddate{#1}}
\def\asciititle#1{\def\theasciititle{#1}}

\def\asciikeywords#1{\def\theasciikeywords{#1}}

\def\shorttitle#1{\def\theshorttitle{#1}}

\let\\\par\let\thelognumber\relax
\let\thevolumenumber\relax\let\thepapernumber\relax
\let\thevolumeyear\relax\let\thesamplenumber\relax\let\startpage\relax
\let\finishpage\relax\let\publishdate\relax\let\receiveddate\relax
\let\reviseddate\relax\let\accepteddate\relax\let\theasciititle\relax
\let\theasciiauthors\relax
\let\theasciiabstract\relax\let\theasciikeywords\relax
\let\theasciiemail\relax\let\theshortauthors\relax\let\theshorttitle\relax

\long\def\maketitlep{   

\count0=\startpage

\gt\hfill      
\beginpicture
\setcoordinatesystem units <0.33truein, 0.33truein> point at 2.2 0.9
\setplotsymbol ({$\cal G$})
\plotsymbolspacing=9truept
\circulararc 315 degrees from 0 1 center at 0 0
\setplotsymbol ({$\cal T$})
\circulararc 315 degrees from 1 -1 center at 1 0
\endpicture
\break
{\small\ifx\thesamplenumber\relax 
Volume \else Sample
\fi\thevolumenumber\ (\thevolumeyear)
\startpage--\finishpage\nl
Published: \publishdate}
\vglue 0.5truein plus 0.4fil minus 0.1truein

{\parskip=0pt\leftskip 0pt plus 1fil\def\\{\par\smallskip}{\ifplaintex\large
\else\Large\fi\bf\thetitle}\par\medskip}   

\vglue 0pt plus 0.1fil 

{\parskip=0pt\leftskip 0pt plus 1fil\def\\{\par}{\sc\theauthors}
\par\medskip}

\vglue 0pt plus 0.1fil 

{\small\parskip=0pt\let\newline\\
{\leftskip 0pt plus 1fil\def\\{\par}{\sl\theaddress}\par}
\expandafter\ifx\theemail\relax    
\relax\else\vglue 5pt plus 0.02fil minus 2pt\def\\{\stdspace{\rm 
and}\stdspace} 
\cl{Email:\stdspace\tt\theemail}\fi
\ifx\theurl\relax                  
\relax\else\vglue 5pt plus 0.02fil minus 2pt\def\\{\stdspace{\rm 
and}\stdspace}
\cl{URL:\stdspace\tt\theurl}\fi\par}

\vglue 7pt plus 0.3fil minus 3pt

{\bf Abstract}
\vglue 5pt plus 0.1fil minus 2pt

\theabstract

\vglue 7pt plus 0.3fil minus 3pt

{\bf AMS Classification numbers}\quad Primary:\quad \theprimaryclass

Secondary:\quad \thesecondaryclass

\vglue 5pt plus 0.3fil minus 2pt

{\bf Keywords:}\quad \thekeywords

\vglue 10pt plus 0.5fil minus 5pt

{\small  Proposed: \theproposer\hfill Received: \receiveddate\nl
Seconded: \theseconders\hfill 
\ifx\reviseddate\relax                         
Accepted: \accepteddate                        
\else
Revised: \reviseddate                          
\fi}
\eject
}       

\let\maketitlepage\maketitlep
\let\maketitle\maketitlepage

\font\phead=cmsl9 scaled 950
\font=cmsl9 scaled 1050
\font\pnum=cmbx10 scaled 913
\font\lnum=cmbx10 
\font\pfoot=cmsl9 scaled 950
\font=cmsl9 scaled 1050
\ifplaintex
\headline{\vbox to 0pt{\vskip -4.5mm\line{\small\phead\ifnum
\count0=\startpage ISSN 1364-0380 (on line)
1465-3060 (printed) \hfill {\pnum\folio}\else\ifodd\count0\def\\{ }%
\ifx\theshorttitle\relax\thetitle\else\theshorttitle\fi\hfill{\pnum\folio}
\else\def\\{ and }{\pnum\folio}\hfill\ifx\theshortauthors\relax\theauthors
\else\theshortauthors\fi\fi\fi}\vss}}
\footline{\vbox to 0pt{\vglue 0mm\line{\small\pfoot\ifnum\count0=\startpage
\copyright\ \gtp\hfill\else
\gt, Volume \thevolumenumber\ (\thevolumeyear)\hfill\fi}\vss
}}
\else
\makeatletter
\def\@oddhead{{\small\lhead\ifnum\count0=\startpage ISSN 1364-0380 (on line)
1465-3060 (printed) \hfill {\lnum\number\count0}\else\ifodd\count0
\def\\{ }\ifx\theshorttitle\relax \thetitle \else\theshorttitle\fi\hfill
{\lnum\number\count0}\else\def\\{ and }{\lnum\number\count0}
\hfill\ifx\theshortauthors\relax 
\theauthors\else\theshortauthors\fi\fi\fi}}\def\@evenhead{\@oddhead}
\def\@oddfoot{\small\lfoot\ifnum\count0=\startpage\copyright\ \gtp\hfill\else
\gt, Volume \thevolumenumber\ (\thevolumeyear)\hfill\fi}
\def\@evenfoot{\@oddfoot}
\makeatother
\fi

\newwrite\gtoutfile
\long\gdef\makeheadfile{  
{\def\\{, }\def\s{ }
\immediate\openout\gtoutfile head.xxx
\immediate\write\gtoutfile{To: math@arxiv.org}
\immediate\write\gtoutfile{Subject: put or rep NNNNN:pppp}
\immediate\write\gtoutfile{--text follows this line--}
\immediate\write\gtoutfile{Proxy-for: \ifx\theasciiauthors\relax
\theauthors\else\theasciiauthors\fi\s<\ifx\theasciiemail\relax\theemail\else\theasciiemail\fi>}
\immediate\write\gtoutfile{\noexpand\\}
\immediate\write\gtoutfile{Authors: \ifx\theasciiauthors\relax
\theauthors\else\theasciiauthors\fi}
\immediate\write\gtoutfile{Title: \ifx\theasciititle\relax
\thetitle\else\theasciititle\fi}
\immediate\write\gtoutfile{Subj-class: GT or SG or MG etc}
\immediate\write\gtoutfile{MSC-class: \theprimaryclass\ifx\thesecondaryclass\relax\else, \thesecondaryclass\fi}
\immediate\write\gtoutfile{Journal-ref: Geom. Topol. \thevolumenumber
(\thevolumeyear) \startpage-\finishpage}
\immediate\write\gtoutfile{Comments: Published by Geometry and Topology at}
\immediate\write\gtoutfile{\s\s http://www.maths.warwick.ac.uk/gt/GTVol\thevolumenumber/paper\thepapernumber.abs.html}
\immediate\write\gtoutfile{\noexpand\\}
\immediate\write\gtoutfile{}
\ifx\theasciiabstract\relax
\immediate\write\gtoutfile{\theabstract}\else
\immediate\write\gtoutfile{\theasciiabstract}\fi
\immediate\write\gtoutfile{}
\immediate\write\gtoutfile{\noexpand\\}
\immediate\write\gtoutfile{}
\immediate\closeout\gtoutfile}}  

\def\maketitlepage{\maketitlep\makeheadfile}
\let\maketitle\maketitlepage

\lognumber{116}
\volumenumber{4}\papernumber{18}\volumeyear{2000}
\pagenumbers{517}{535}
\proposed{Dieter Kotschick}
\seconded{Robion Kirby, Yasha Eliashberg}
\received{12 April 2000}
\revised{8 December 2000}
\accepted{17 December 2000}
\published{21 December 2000}

\usepackage{amssymb,amsmath,ifthen}

\makeatletter
\newcommand{\opname}[2][it]%
        {%
	\let\iflong@one=\iftrue
	\ifthenelse{\equal{#1}{it}}
	   {\expandafter\def\csname #2\endcsname{\ensuremath{\mathit{#2}}}
	    \let\iflong@one=\iffalse}{}
	\ifthenelse{\equal{#1}{bb}}
	   {\expandafter\def\csname #2\endcsname{\ensuremath{\mathbb{#2}}}
	    \let\iflong@one=\iffalse}{}
	\ifthenelse{\equal{#1}{rm}}
	   {\expandafter\def\csname #2\endcsname{\ensuremath{\mathrm{#2}}}
	    \let\iflong@one=\iffalse}{}
	\ifthenelse{\equal{#1}{bf}}
	   {\expandafter\def\csname #2\endcsname{\ensuremath{\mathbf{#2}}}
	    \let\iflong@one=\iffalse}{}
	\ifthenelse{\equal{#1}{tt}}
	   {\expandafter\def\csname #2\endcsname{\ensuremath{\mathtt{#2}}}
	    \let\iflong@one=\iffalse}{}
	\ifthenelse{\equal{#1}{cal}}
	   {\expandafter\def\csname #2\endcsname{\ensuremath{\mathcal{#2}}}
	    \let\iflong@one=\iffalse}{}
	\ifthenelse{\equal{#1}{bit}}
	   {\expandafter\def\csname #2\endcsname{\ensuremath{\text{\boldmath $#2$}}}
	    \let\iflong@one=\iffalse}{}
	\iflong@one 
	   \expandafter\def\csname #2\endcsname{\ensuremath{{#1}}}
	\fi
	}
\makeatother

\opname[rm]{genus}
\opname[rm]{Center}
\opname[rm]{rk}
\opname[it]{pr}
\opname[rm]{ker}
\opname[it]{Mon}
\opname[it]{Map}
\opname[it]{Id}
\opname[it]{Im}
\opname[rm]{Out}
\opname[bb]{Z}
\opname[bb]{R}
\opname[bb]{C}
\opname[bb]{P}
\opname[cal]{M} 
\opname[cal]{G}
\opname[cal]{S}
\opname[cal]{P}
\opname[cal]{B}
\opname[\C P]{CP}
\opname[it]{Diff}
\opname[\mathbf{Diff}]{Diffbf}
\opname[\mathbf{Map}]{Mapbf}
\let\text=\mbox

\newcommand{\larrow}[1]{\stackrel{#1}{\rightarrow}}
\newcommand{\pgn}{P_{g,n}}
\newcommand{\pf}{{\pi_1(F_{b_0},x_0)}}
\newcommand{\pim}{{\pi_1(M^3)}}
\newcommand{\pb}{{\pi_1(B,b_0)}}
\newcommand{\px}{{\pi_1(X,x_0)}}

\let\diffeo=\cong

\let\phi=\varphi

\newtheorem{thm}{Theorem}[section]

\newtheorem{lem}{Lemma}[section]
\newtheorem{prop}{Proposition}[section]
\newtheorem*{claim}{Claim}

\theoremstyle{definition}

\newtheorem*{definition}{Definition}

\begin{document}
\title{Symplectic Lefschetz fibrations on $S^1\times M^3$}
\shorttitle{Symplectic Lefschetz fibrations on a circle times a 3-manifold}
\asciititle{Symplectic Lefschetz fibrations on S^1 x M^3}

\author{Weimin Chen\\Rostislav Matveyev}

\address{UW-Madison, Madison, WI 53706, 
USA\\SUNY at Stony Brook, NY 11794, USA}
\email{wechen@math.wisc.edu\\slava@math.sunysb.edu}

\begin{abstract}
In this paper we classify symplectic Lefschetz fibrations (with empty
base locus) on a four-manifold which is the product of a three-manifold
with a circle. This result provides further evidence in support of the
following conjecture regarding symplectic structures on such a
four-manifold: if the product of a three-manifold with a circle admits a
symplectic structure, then the three-manifold must fiber over a
circle, and up to a self-diffeomorphism of the four-manifold, the
symplectic structure is deformation equivalent to the canonical
symplectic structure determined by the fibration of the three-manifold
over the circle.
\end{abstract}

\asciikeywords{Four-manifold, symplectic structure, Lefschetz fibration, 
Seiberg-Witten invariants}



\maketitlepage

\section{Introduction and statement of results}

Suppose $M^3$ is a closed oriented 3--manifold.  If $M^3$ fibers over
$S^1$, then the 4--manifold $X=S^1\times M^3$ has a symplectic
structure canonical up to deformation equivalence.

An interesting question motivated by Taubes' fundamental research on
symplectic 4--manifolds
\cite{T} asks whether the converse is true:
{\em
Does every symplectic structure (up to deformation equivalence)
on $X=S^1\times M^3$
come from a fibration of $M^3$ over $S^1$,
in particular, is it true that $M^3$ must be fibered?
}

In this paper we prove the following:

\begin{thm} Let $M^3$ be a closed 3--manifold which contains no fake
3--cell.
If $\omega$ is a symplectic structure on the 4--manifold
$X=S^1\times M^3$ defined through a Lefschetz fibration
(with empty base locus),
then:
\begin{enumerate}
\item $M^3$ is fibered, with a fibration $p\co M^3\rightarrow S^1$.
\item There is a self-diffeomorphism $h\co X\rightarrow X$
such that $h^*\omega$ is deformation equivalent to the 
canonical symplectic structure on $X$ associated to
the fibration $p\co M^3\rightarrow S^1$.
\end{enumerate}
\end{thm}

{\bf Remarks}\qua
(1)\qua A 3--manifold could admit essentially different fibrations
over a circle (see eg \cite{McT}). The fibration 
$p\co M^3\rightarrow S^1$ in the theorem is specified by the Lefschetz 
fibration.

(2)\qua The self-diffeomorphism $h\co X\rightarrow X$ 
is homotopic to the identity.

(3)\qua R Gompf and A Stipsicz have shown \cite{GS} that 
for any 4--manifold with the rational homology of $S^1 \times S^3$, any
Lefschetz pencil or fibration (even allowing singularities with the
wrong orientation) must be a locally trivial torus fibration over $S^2$, in
particular, the manifold is $S^1 \times L(p,1)$ with the obvious fibration.

\medskip

A stronger version of the theorem, in which we also classify 
symplectic Lefschetz fibrations on $S^1\times M^3$, is given in section 4.

The proof of the theorem consists of two major steps.  First, we show
that any symplectic Lefschetz fibration on $X = S^1 \times M^3$ has no
singular fibers, ie, it is a locally trivial fibration. This result
is the content of lemma 3.4 (Lemma on Vanishing Cycles). Secondly, we
show that every such fibration of $X = S^1 \times M^3$ is induced in a
certain way from a fibration of $M^3$ over a circle.

The essential ingredient in the proof of the Lemma on Vanishing Cycles is
a result of D Gabai, which says roughly that the minimal genus of
an immersed surface representing a given homology class in a
3--manifold is equal to the minimal genus of an embedded surface
representing the same homology class. This result is specific for
dimension three and does not hold in dimension four.

This paper is organized as follows: Section 2 contains a brief review
of some background results.  Section 3 is devoted to a preparatory
material for the proof of the main theorem, in particular, it contains
the Lemma on Vanishing Cycles.  In section 4, the full version of the main
theorem is stated and proved.

\medskip

{\bf Acknowledgments}\qua The authors are grateful to John McCarthy
for pointing out an erroneous quotation of the Nielsen Representation
Theorem in an early version.  The second author is indebted to
Siddhartha Gadgil for numerous discussions on 3--dimensional topology.
The authors wish to thank the editors and referees for their comments
which have helped to improve the presentation greatly.  The first
author is partially supported by an NSF grant.

\section{Recollections}

\subsection{Lefschetz pencils and fibrations}

Let $X$ be a closed, oriented, smooth 4--manifold. A Lefschetz pencil
on $X$ is a smooth map $P\co X\setminus \B\rightarrow \CP^1$ defined
on the complement of a finite subset $\B$ of $X$, called the base
locus, such that each point in $\B$ has an orientation-preserving
coordinate chart in which $P$ is given by the projectivization map
$\C^2\setminus\{0\}\rightarrow \CP^1$, and each critical point has an
orientation-preserving chart on which $P(z_1,z_2)=z_1^2+z_2^2$.
Blowing up at each point of $\B$, we obtain a Lefschetz fibration on
$X\# n\overline{\CP^2}$ ($n=\# \B$) over $\CP^1$ with fiber
$F_t=P^{-1}(t)\cup \B$ for each $t\in \CP^1$.

More generally, a Lefschetz fibration on a closed oriented smooth 4--manifold
$X$ is a smooth map $P\co  X\rightarrow B$ where $B$ is a Riemann
surface, such that each critical point of $P$ 
has an orientation-preserving chart on which $P(z_1,z_2)=z_1^2+z_2^2$.
We require that each fiber is connected and contains at most one critical 
point.
Every Lefschetz fibration can be changed to satisfy these two conditions.
A Lefschetz fibration is called symplectic if there is a symplectic structure
$\omega$ on $X$ whose restriction to each regular fiber is non-degenerate.

In an orientation-preserving chart at a critical point $x\in X$, 
the map $P$ is given by $P(z_1,z_2)=z_1^2+z_2^2$.
Let $t\in\R\subset\C$ be a positive
regular value.
The fiber $F_t=P^{-1}(t)$ contains a simple closed loop $\gamma$ which is
the intersection of $F_t$ with the real plane $\R^2\subset \C^2$, ie, the
boundary of the disc in $X$ defined by $x_1^2+x_2^2\leq t$.
This
simple closed loop $\gamma$ on $F_t$ is called the vanishing cycle associated
to the critical point $x$.
A regular neighborhood of the
singular fiber $F_0$ can be described as the 
result of attaching a 2--handle along
the vanishing cycle $\gamma$ to a regular neighborhood of the regular fiber
$F_t$ with framing $-1$ relative to the product framing on $F_t\times S^1$.
The monodromy around the critical value $P(x)$ is a right Dehn twist along the 
vanishing cycle $\gamma$.

We quote an observation of Gompf which roughly says that most of the
Lefschetz fibrations are symplectic.
For a proof, see \cite{GS} or \cite{ABKP}.

\begin{thm}[Gompf]
A Lefschetz fibration on an oriented 4--manifold $X$ is symplectic if the 
fiber class is non-zero in $H_2(X;\R)$.
Any Lefschetz pencil (with non-empty base locus) is symplectic.
\end{thm}

A remarkable theorem of Donaldson \cite{D} says
that any symplectic 4--manifold admits a Lefschetz pencil by symplectic
surfaces.

The following lemma is known for general Lefschetz fibrations. For
completeness, we give a short proof for the case of non-singular
fibrations. This is sufficient for our purpose.

\begin{lem}
Let $F\hookrightarrow X \larrow{P} B$ be a fibration and $\omega_1$,
$\omega_2$ symplectic forms on $X$ with respect to which each fiber
of $P$ is symplectic.
Then $\omega_1$ and  $\omega_2$ are deformation
equivalent if they induce the same orientation on the fiber.
\end{lem}

\proof
We orient the base $B$ so that the pull-back of the volume form $\omega_B$
of $B$, $P^\ast\omega_B$, has the property that 
$\omega_1\wedge P^\ast\omega_B$ is positive on $X$ with respect to 
the orientation induced by $\omega_1$.
Then
$$
\omega(s):=\omega_1+sP^\ast\omega_B
$$
is a symplectic form on $X$ for any $s\geq 0$.
Let $\omega(s,t)=t\omega_2+
(1-t)\omega(s)$ for $0\leq t\leq 1$, $s\geq 0$. Then
\begin{eqnarray*}
\omega(s,t)\wedge\omega(s,t) & = & t^2\omega_2\wedge \omega_2 + 
                                  (1-t)^2\omega(s)\wedge\omega(s)\\
                             &   & +2t(1-t)\omega_2\wedge\omega(s)
\end{eqnarray*}
is positive for all $0\leq t\leq 1$ when $s$ is sufficiently large, since
$\omega_2\wedge\omega(s)$ is positive for large enough $s$.
This implies that $\omega_1$ and $\omega_2$ are deformation equivalent.
\endproof

\subsection{Seiberg--Witten theory}

On an oriented Riemannian 4--manifold $X$, a $Spin^c$ structure $\S$ consists
of a hermitian vector bundle $W$ of rank $4$, together with a Clifford
multiplication $\rho\co T^\ast X\rightarrow End(W)$.
The bundle $W$ decomposes
into two bundles of rank $2$, $W^+\oplus W^-$, such that $\det W^+=\det W^-$.
Here $W^-$ is characterized as the subspace annihilated by $\rho(\eta)$
for all self-dual 2--forms $\eta$.
We write $c_1(\S)$ for the first Chern
class of $W^+$.
The Levi--Civita connection on $X$ coupled
with a $U(1)$ connection $A$ on $\det W^+$ defines a Dirac operator 
$ D_A\co  \Gamma(W^+)\rightarrow \Gamma(W^-)$ from the space of 
smooth sections of $W^+$ into that of $W^-$.

The 4--dimensional Seiberg--Witten equations are the following pair of equations
for a section $\psi$ of $W^+$ and a $U(1)$ connection $A$ on $\det W^+$:
$$
\begin{array}{ccc}
\rho(F_A^+)-\{\psi\otimes\psi^\ast\} & = & 0\\
D_A\psi                              & = & 0\\
\end{array}
$$
Here $F^+$ is the projection of the curvature onto the self-dual forms, and
the curly brackets denote the trace-free part of an endomorphism of $W^+$.

The moduli space $\M_{\S}$ is the space of solutions $(A,\psi)$ modulo the 
action of the gauge group $\G=Map(X,S^1)$, which is compact with virtual
dimension 
$$
d({\S})=\frac{1}{4}(c_1({\S})^2[X]-2\chi(X)-3\sigma(X)),
$$
where $\chi(X)$ and $\sigma(X)$ are the Euler characteristic and signature 
of $X$ respectively.
When $b^+(X)\geq 1$, for a generic perturbation of the Seiberg--Witten
equations where $\eta$ is a self-dual 2--form
$$
\begin{array}{ccc}
\rho(F_A^++i\eta)-\{\psi\otimes\psi^\ast\} & = & 0\\
D_A\psi                              & = & 0,\\
\end{array}
$$
the moduli space $\M_{{\S},\eta}$ is a compact, canonically oriented, smooth
manifold of dimension $d({\S})$, which contains no reducible solutions
(ie, solutions with $\psi\equiv 0$).
The fundamental class of $\M_{{\S},\eta}$
evaluated against some universal characteristic classes defines the 
Seiberg--Witten invariant $SW({\S})\in\Z$, which is independent of the 
Riemannian metric and the perturbation $\eta$ when $b^+(X)>1$.
When $b^+(X)=1$,
$SW({\S})$ is well-defined if $c_1({\S})^2[X]\geq 0$ and $c_1({\S})$ is not 
torsion, by choosing $||\eta||$ sufficiently small.
The set of complex line
bundles $\{E\}$ on $X$ acts on the set of $Spin^c$ structures freely and 
transitively by $(E,{\S})\rightarrow {\S}\otimes E$.
We will call a cohomology class $\alpha\in H^2(X;\Z)$
a Seiberg--Witten basic class if there exists a $Spin^c$ structure $\S$ such
that $\alpha=c_1({\S})$ and $SW({\S})\neq 0$.
There is an involution $I$
acting on the set of $Spin^c$ structures on $X$ which has the property that
$c_1({\S})=-c_1(I({\S}))$ and $SW({\S})=\pm SW(I({\S}))$.
As a consequence,
if a cohomology class $\alpha\in H^2(X;\Z)$ is a Seiberg--Witten basic class,
so is $-\alpha$.

We will use the following fundamental result:

\begin{thm}[Taubes]
Let $(X,\omega)$ be a symplectic 4--manifold with canonical line
bundle $K_\omega$.
Then $c_1(K_\omega)$ is a Seiberg--Witten basic class if $b^+(X)>1$ or
$b^+(X)=1$, $c_1(K_\omega)\cdot [\omega]>0$ and $2\chi(X)+3\sigma(X)\geq 0$.
\end{thm}

\subsection{Gabai's Theorem}

The following theorem of Gabai says that, given a singular oriented
surface in a closed oriented 3--manifold, one can find an embedded surface (not
necessarily connected) representing the same homology class and having
the same topological complexity as the singular surface.

Let us recall the definition of Thurston norm and the singular norm on
the second homology group of a compact 3--manifold.
Let $S$ be an orientable surface. The complexity of $S$ is defined by
$x(S)=\sum_{S_i}\max(-\chi(S_i), 0)$, where the summation is taken over
connected components of $S$.
For a closed, oriented 3--manifold $M$, the Thurston norm $x(z)$ and
singular norm $x_s(z)$ of a 
homology class $z\in H_2(M;\Z)$ are defined by 
\begin{eqnarray*}
x(z) &=&\min\{x(S)|\text{$S$ is an
embedded surface representing $z$}\},\\
x_s(z)&=&\inf\{\frac{1}{n} x(S)|f\co S\rightarrow M, f([S])=nz\}.
\end{eqnarray*}

\begin{thm}[Gabai \cite{Ga}]
Let $M$ be a closed oriented 3--manifold. Then
$x_s(z)=x(z)$ for all $z\in H_2(M;\Z)$.
\end{thm}

Gabai's theorem is specific for dimension three and fails in dimension
four in general. Precisely because our 4--manifold under consideration
is the product of a 3--manifold with a circle, we are able to apply
Gabai's theorem to yield a stronger estimate for the 4--manifold. This
is the essential point in the proof of the Lemma on Vanishing Cycles.

\section{Preparatory material for the proof of the theorem}

This section is devoted to preparatory 
material necessary for the proof of 
the theorem.

\subsection{$\Mapbf(F)$ and $\Diffbf(F)$}

We first list some results about the mapping class group $\Map(F)$
and the diffeomorphism group $\Diff(F)$ of an orientable surface $F$.

The reader is referred to \cite{Ker} for the proof of the following
proposition.

\begin{prop}{\rm(Nielsen Representation Theorem)}\ \
Every finite subgroup $G$ of the mapping 
class group of a big ($\geq 2$)
genus surface can be lifted to the 
diffeomorphism group of that surface.
Moreover, there is a conformal structure on $F$ such that 
$G$ is realized by conformal isometries of $F$.
\end{prop}

We will also need an analogue of the Nielsen Representation Theorem in a 
slightly different setting. It states that a commutativity relation between 
an arbitrary element and a finite order element of the mapping class group 
of a big ($\geq 2$) genus surface can be lifted to the homeomorphism group 
of the surface.

\begin{lem}
Let $F$ be a surface with $\genus(F)\geq 2$ and $\psi,\phi$ be two
mapping classes of $F$ such that $\phi$ has finite order and 
$\psi\circ\phi=\phi\circ\psi$ in $\Map(F)$.
Then there are homeomorphisms $\Psi,\Phi\co F\rightarrow F$ in the mapping 
classes $\psi,\phi$ respectively, such that $\Phi$ has finite order, 
and $\Psi\circ\Phi=\Phi\circ\Psi$ as homeomorphisms.
\end{lem}

\proof The proof is based on a theorem of 
Teichm\"uller \cite{Tm1,Tm2}. A very clear treatment of Teichm\"uller's
theorem can be found in \cite{B}.

We first recall the notions of quasi-conformal mappings and total dilatation
of homeomorphisms of a surface.
Let $F$ be a Riemann surface with a given conformal structure,
and $f\co F\rightarrow F$ be an orientation 
preserving homeomorphism. At a point $p\in F$ where $f$ is $C^1$--smooth 
we can measure the deviation of $f$ from being conformal by the ratio of 
the bigger axis to the smaller axis of an infinitesimal ellipse, 
which is the image of an infinitesimal circle around the point $p$ under $f$.
This ratio is called the local dilatation of $f$ and is denoted by $K_p[f]$.
A homeomorphism $f\co F\rightarrow F$ is called a quasi-conformal mapping
if $K_p[f]$ is defined for almost all $p\in F$ and $\sup K_p[f]< \infty$.
The number $K[f]=\sup K_p[f]$ is called the total dilatation of the 
quasi-conformal mapping $f$.

Teichm\"uller's theorem states that {\sl among all the homeomorphisms
in a given mapping class of a Riemann surface $F$ of $\genus(F)\geq 2$,
there is a unique one which minimizes the total dilatation.}

Now back to the proof of the lemma. By the Nielsen Representation Theorem
(cf Proposition 3.1), the mapping class $\phi$ can be represented by an
isometry $\Phi$ of $F$ with respect to some conformal structure $\alpha$.
By Teichm\"uller's theorem, there is a unique homeomorphism $\Psi$ 
in the mapping class $\psi$ which minimizes the total dilatation with respect 
to the conformal structure $\alpha$. 
Since $\Phi$ is an isometry of $\alpha$, $\Phi^{-1}\circ\Psi\circ\Phi$ has
the same total dilatation as $\Psi$. On the other hand, 
$\Phi^{-1}\circ\Psi\circ\Phi$ and $\Psi$ are in the same mapping class since
$\psi\circ\phi=\phi\circ\psi$ in $\Map(F)$. 
The uniqueness of the extremal homeomorphism in Teichm\"uller's theorem 
implies that they must coincide, and $\Psi\circ\Phi=\Phi\circ\Psi$ 
as homeomorphisms.
\endproof

\begin{prop} Let $F_g$ be a closed orientable 
surface of genus $g$, and denote by $\Diff_0(F_g)$ the identity component of
$\Diff(F_g)$.
Then for all $g\geq 2$,
$\Diff_0(F_g)$ is contractible,
for $g=1$, 
$\Diff_0(F_g)$ is homotopy equivalent to the identity component of the group
of conformal automorphisms of the torus and, 
for $g=0$, $\Diff_0(F_g)$ is homotopy
equivalent to $SO(3)$.
\end{prop}

The reader is referred to \cite{EE} for this result.

The rest of this subsection concerns locally trivial fibrations of
4--manifolds over a Riemann surface. 

\begin{definition} 
Two locally trivial fibrations $F\hookrightarrow X \larrow{P} B$ and 
$F'\hookrightarrow X^\prime \larrow{P'} B'$
are said to be equivalent if there are diffeomorphisms 
$f\co B\rightarrow B'$ and 
$\tilde{f}\co  X\rightarrow X^\prime$, such that the following diagram
$$ 
\begin{array}{rcl}
X&\larrow{\tilde{f}}&X^\prime\\
{\scriptstyle P}\downarrow\;&&\:\downarrow {\scriptstyle P'}\\
B&\larrow{f}&B'\\
\end{array}
$$
commutes.
When $X=X^\prime$, we say that $P$ and $P^\prime$ are 
strongly equivalent if $\tilde{f}$ is isotopic to the identity.
\end{definition}

Each fibration $P$ defines a monodromy homomorphism $\Mon_P$ 
from the fundamental group of the base 
to the mapping class group of the fiber.

\begin{lem}
When the fiber is high-genus ($\geq 2$), monodromy $\Mon_P$ determines the
equivalence class of a fibration $P$.
\end{lem}

\proof 
We write the base $B$ as $B_0\cup D$ where $B_0$ is
a disc with several 1--handles attached and $D$ is a disc.
We choose a base point $b_0\in B_0\cap D$.
The fibration $P$
restricted to $B_0$ is determined by $\Mon_P$ as a representation of 
$\pi_1(B_0,b_0)$ in the mapping class group of the fiber $F_{b_0}$ over
the base point $b_0$.
Over the
disc $D$, $P$ is trivial.
We can recover $P$ on 
the whole manifold $X$ by gluing $P^{-1}(B_0)$ and $P^{-1}(D)$ along the 
boundary via some gluing map $\phi$, and $P$ depends only on the homotopy
class of $\phi$ viewed as a map $S^1\rightarrow \Diff_0(F_{b_0})$.
Now we
recall the fact that $\Diff_0(F_g)$ is 
contractible for all $g\geq 2$, so that $\phi\sim \Id$ when $\genus(F_{b_0})
\geq 2$.
The lemma follows.
\endproof

\subsection{Cyclic coverings and $S^1$--valued functions}

Here we describe a construction, which,  starting with a $S^1$--valued
function on a $CW$--complex and a positive integer,  gives
a cyclic finite covering of the $CW$--complex 
and a generator of the deck transformation group.
An inverse of this construction is then discussed.

For the rest of the paper we view $S^1$ as the unit circle in $\C$,
oriented in the usual way, ie, in the counter-clockwise direction.

Consider a $CW$--complex $Y$ and let
 $g\co Y\rightarrow S^1$ be a function.
It determines a cohomology class $[g]\in H^1(Y;\Z)$.
Denote by $d(g)$ the divisibility of $[g]$.
For a positive integer $n$ and a function $g\co Y\rightarrow S^1$ 
such that $\gcd(n,d(g)) = 1$, define a subset 
$\widetilde{Y}\subset S^1\times Y$ by 
$$
\widetilde{Y}=\{(t, y)\in S^1\times Y\;|\; t^n=g(y)\}.
$$
Obviously $\widetilde{Y}$ is invariant under the transformation 
$\phi\co S^1\times Y\rightarrow S^1\times Y$ induced
from rotation of $S^1$ by angle 
$\frac{2\pi}{n}$ in the positive direction.
Let $pr\co \widetilde{Y}\rightarrow Y$ be the restriction of the projection 
from $S^1\times Y$ onto the second factor.
It is a cyclic $n$--fold covering of $Y$
and $\phi$ generates the group of deck transformations
of this covering.

Now we would like to invert this construction, ie, starting 
with a finite cyclic covering $pr\co \widetilde{Y}\rightarrow Y$ and a generator
$\phi\co \widetilde{Y}\rightarrow\widetilde{Y}$ 
of the structure group of $pr$,
find a function $g\co \widetilde{Y}\rightarrow S^1$
such that the construction above yields  $pr\co \widetilde{Y}\rightarrow Y$ and 
$\phi$.
However this is not always possible.

Fix a finite cyclic covering $\pr\co \widetilde{Y}\rightarrow {Y}$ with 
a generator $\phi\co \widetilde{Y}\rightarrow\widetilde{Y}$ 
of the structure group of $pr$.
Let $G=\langle\phi\rangle$ be the group of deck transformations of $pr$ 
and $n = |G|$.
Define an action of $G$
on $S^1$ by $\phi\cdot z = e^{2\pi i/n}z$.
This gives rise to an $S^1$--bundle $Z$ over $Y$
and $\widetilde{Y}$ sits naturally in $Z$.
The Euler class of this bundle is a torsion class, for the 
$n^\mathrm{th}$ power
of the bundle (considered as $U(1)$ bundle) is trivial.
Moreover, there is a choice of 
trivialization, canonical up to rotations of $S^1$, given by an 
$n$--valued section of $pr$, which becomes a section
in the $n^\mathrm{th}$ power of the bundle.
Suppose the Euler class vanishes, then $Z\diffeo S^1\times Y$
and this diffeomorphism is canonical up to isotopy.
Consider the map $\tilde{g}\co \widetilde{Y}\rightarrow S^1$,
which is the restriction to $\widetilde{Y}$ 
of the projection from $S^1\times Y$ onto the $S^1$--factor.
The map $\tilde{g}^n$ is $\phi$ invariant and therefore 
descends to a map $g\co Y\rightarrow S^1$.
It is an easy exercise to show that the construction above applied to the
pair $(g,n)$ yields $(pr, \phi)$.

Thus we have proved the following:

\begin{lem}
Let $Y$ be a $CW$--complex such that $H^2(Y;\Z)$ has no torsion.  Then the
construction above gives a 1--1 correspondence between two sets
$A=\{(pr,\:\phi)\}$ and $B=\{(n,[g])\}$,
where $pr\co \widetilde{Y}\rightarrow Y$ is a finite cyclic covering, 
$\phi\co \widetilde{Y}\rightarrow \widetilde{Y}$ is a generator of the deck 
transformation group of $pr$, $n$ is a positive integer and $[g]$ is the 
homotopy class of a map $g\co Y\rightarrow S^1$ with $\gcd(n,d(g)) = 1$.
\end{lem}

\subsection{Surface fibrations over a torus}

In this subsection we will explore a family of surface fibrations
of $X=S^1\times M^3$ over a torus, where $M^3$ is a closed orientable 
3--manifold fibered over a circle with fiber $\Sigma$ and fibration
$p\co M^3\rightarrow S^1$.
Denote by $\Sigma_\beta$ the fiber over a point 
$\beta\in S^1$ and by $[p]\in H^1(M^3;\Z)$ the homotopy class of $p$.

Consider a smooth function $g\co M^3\rightarrow S^1$ and
denote by $[g]\in H^1(M;\Z)$ its homotopy class
and by $d(g)$ the divisibility of $[g]$ restricted to the fiber $\Sigma$.
Let $n$ be a positive integer such that $\gcd(n, d(g)) = 1$.

Define $P_{g,n}\co X=S^1\times M^3\rightarrow S^1\times S^1$ by
$P_{g,n}(t, m) = ( t^n\overline{g(m)},p(m))$, where
$\overline{g(m)}$ stands for the complex conjugate
of $g(m)$ in $S^1$.
The map $\pgn$ is a locally trivial fibration of $X$ over a torus
such that the fiber $F_{(\alpha,\beta)}$ over a point 
$(\alpha, \beta)\in S^1\times S^1$ 
is the graph of a multi-valued function
$(\alpha g|_{\Sigma_\beta})^{\frac{1}{n}}$ on $\Sigma_\beta$, ie, 
$$
F_{(\alpha,\beta)}=
\{(t,m)\in S^1\times \Sigma_\beta\: |\: t^n=\alpha g(m)\}\subset 
S^1\times\Sigma_\beta\subset S^1\times M^3.
$$

Fix $g\co M^3\rightarrow S^1$ and $n$ as above.
Let us find the monodromy of $\pgn$.
First of all, the projection $S^1\times\Sigma_\beta\rightarrow \Sigma_\beta$
restricted to $F_{(\alpha,\beta)}$ is a cyclic n--fold covering
$\pr\co F_{(\alpha,\beta)}\rightarrow\Sigma_\beta$.
Denote by $\phi$ the self-diffeomorphism of $F_{(1,1)}$ induced by
a rotation of the $S^1$--factor in $S^1\times\Sigma_1$ by angle $2\pi/n$,
then $\phi$ generates the group of deck transformations of 
$\pr\co F_{(1,1)}\rightarrow\Sigma_1$.
Secondly,
denote by $\Mon_p\in \Map(\Sigma_1)$ the monodromy of
$p\co M^3\rightarrow S^1$, then $\Mon_p$ pulls back to an element
in $Map(F_{(1,1)})$.

Let $q=S^1\times\{1\}$ and $r=\{1\}\times S^1$ be the ``coordinate'' simple 
closed loops in $S^1\times S^1$.
Then it is easily seen that 
the monodromy $\Mon_{\pgn}(q)$ of $\pgn$ along $q$ is equal to 
$[\phi]$ and the monodromy $\Mon_{\pgn}(r)$ along $r$ is the 
pull-back of $\Mon_p$ to $\Map(F_{(1,1)})$.

We end this subsection with the following classification of $\pgn$.

\begin{prop}
Two fibrations $P_{g,n}$ and $P_{g',n'}$ are strongly equivalent if and 
only if $n=n'$ and $[g]\equiv[g']\!\!\mod([p])$, ie,  when 
the two planes in $H^1(X;\Z)=H^1(S^1;\Z)\oplus H^1(M^3;\Z)$
spanned by two pairs of vectors $(n[t] - [g], [p])$ and $(n'[t] - [g'], [p])$
coincide.
\end{prop}

\proof
The necessity follows from homotopy theoretical considerations:
the two said planes are pull-backs of $H^1(S^1\times S^1;\Z)$ by $P_{g,n}$ 
and $P_{g',n'}$ respectively, therefore they
must coincide if $P_{g',n'}$ is homotopic to 
$P_{g,n}$ post-composed with a self-diffeomorphism
of $S^1\times S^1$.

On the other hand, by changing basis on one of the tori, we can
arrange that the homotopy classes of $g$ and $g'$ are equal.
Then a homotopy from $g$ to $g'$ 
leads to a strong equivalence between $P_{g,n}$ and $P_{g',n'}$.
\endproof

We observe that when the fiber $\Sigma$ of $p\co M^3\rightarrow S^1$
has genus zero, there is a unique equivalence class of $P_{g,n}$, ie,
the trivial one when $[g]=0$ and $n=1$.

\subsection{Lemma on Vanishing Cycles}

\begin{lem} Let $M^3$ be a closed 3--manifold and 
$P\co X=S^1\times M^3 \rightarrow B$ be a symplectic Lefschetz fibration
with regular fiber $F$.
Then $P$ has no singular fibers.
\end{lem}

\proof
We first observe that every vanishing cycle
must be non-separating since $X$ has an even intersection form.
There are three cases to consider according to the genus of the fiber.

(1)\qua The fiber $F$ is a sphere. There are no vanishing cycles since every
    curve on a sphere is separating.

(2)\qua The fiber $F$ is a torus. Recall a formula for the Euler characteristic
of the total space of a Lefschetz fibration
$$
\chi(X)=\chi(F)\cdot\chi(B)+\#\text{\{vanishing cycles\}}.
$$
There must be no vanishing cycles since $\chi(F)=0$ and $\chi(X)=0$.

(3)\qua The fiber $F$ is a high genus ($\geq 2$) surface. We first show that
    the canonical class $K$ is Seiberg--Witten basic. By Theorem 2.2,
    we only need to show that, if $b^+(X)=1$, then $K\cdot [\omega]
    >0$, where $\omega$ is the symplectic form on $X$. This could be
    seen as follows. The 4--manifold $X$ has a hyperbolic intersection
    form since $b^+(X)=1$ and $\sigma(X)=0$. Let $x,y\in H^2(X;\R)$
    form a hyperbolic basis, ie, $x\cdot x = y\cdot y =0$ and $x\cdot
    y = 1$. Since $F\cdot F = 0$, $K\cdot K = 2\chi (X) + 3\sigma(X) =
    0$ and $K\cdot F=2g_F-2>0$, we can assume without loss of
    generality that $PD[F] = ax$ and $K=by$ for some positive $a$ and
    $b$. Let $[\omega]=\alpha x+\beta y$. Observe that
    $[\omega]\cdot[\omega] = \alpha\beta >0$ and $[\omega]\cdot F = a\beta
    >0$, therefore $\alpha$ and $\beta$ are both positive. But
    $K\cdot[\omega] = b\alpha>0$. Thus, by Theorem 2.2, $K$ is a
    Seiberg--Witten basic class.

    Let $V\in H_2(M^3;\Z)$ be the homology class of the projection of
    $F$ into $M^3$. Suppose there is a singular fiber. Since the
    vanishing cycle for this singular fiber is non-separating, the
    class $V$ can be represented by a map $f\co F^\prime\rightarrow
    M^3$ such that $g_{F^\prime}=g_F-1$. The singular norm of $V$,
    $x_s(V)$, is less than or equal to $2g_{F^\prime}-2$. By Theorem
    2.3, there are embedded surfaces $S_i$ in $M^3$ such that $\sum_i
    [S_i]=V$ and $\sum_i x(S_i)= x_s(V)\leq 2g_{F^\prime}-2$, where
    $x$ stands for the complexity of an orientable surface.  On the
    other hand, since $K$ is Seiberg--Witten basic, by the adjunction
    inequality (cf \cite{K}), we have
$$
x(S_i) \geq |K\cdot S_i|, \hspace{2mm}\text{for all $i$}.
$$
Therefore,
$$
2g_{F^\prime}-2\geq \sum_i x(S_i)\geq |K\cdot V|=K\cdot F=2g_F-2,
$$
which is a contradiction.  We used the fact that $K\cdot V=K\cdot
F$. This is because $H_1(S^1;\Z)\otimes H_1(M^3;\Z)\subset H_2(X;\Z)$
consists of classes represented by embedded tori in $X$, thus, by the 
adjunction inequality, $K\cdot H=0$ for any $H\in H_1(S^1;\Z)\otimes
H_1(M^3;\Z)$. This concludes the proof for the case of high genus
fiber.
\endproof

\section{The main theorem}

Let us recall that if $M$ is a closed oriented 3--manifold fibered
over $S^1$, then the 4--manifold $X=S^1\times M$ carries a canonical
(up to deformation) symplectic structure compatible with the
orientation on $X$.  Note that we have canonically oriented $S^1$ so
that the orientation of $M$ determines an orientation of the fiber of
$M\rightarrow S^1$.

The Poincar\'e associate of a 3--manifold $M$, denoted by $\P(M)$, is defined
by the condition that $\P({M})$ contains no fake 3--cell and 
$M =\P({M})\# A$ where $A$ is a homotopy 3--sphere. The theorem about unique 
normal prime factorization of 3--manifolds implies that the Poincar\'e 
associate exists and is unique \cite{H}.
An orientation on $M$ canonically determines an orientation on 
the Poincar\'e associate $\P(M)$.

\begin{thm}
Let $M^3$ be a closed oriented 3--manifold and $X=S^1\times M^3$.
Let $P\co X\rightarrow B$ be a symplectic Lefschetz fibration with respect to 
some symplectic form $\omega$ on $X$ compatible with the orientation.
Denote by $F$ the regular fiber of $P$.
Then:
\begin{enumerate}
\item The Poincar\'e associate $\P({M}^3)$ of $M^3$ is
fibered over $S^1$ with fibration $p\co \P({M}^3)\rightarrow S^1$.
\item There is a diffeomorphism $h\co S^1\times \P({M}^3)\rightarrow X$
such that $h^*\omega$ is deformation equivalent to the 
canonical symplectic structure on $S^1\times\P(M^3)$ defined through the 
fibration $p\co \P(M^3)\rightarrow S^1$.
\end{enumerate}

There are two possibilities for the Lefschetz fibration $P$:
\begin{itemize}
\item [\rm(a)] If $P(S^1\times\{pt\})$ is not null-homotopic in
$B$, then $\genus(B) = 1$, and $P\circ h=P_{g,n}$ for some
integer $n>0$ and map $g\co \P(M^3)\rightarrow S^1$
(see subsection 3.3 for the definition of $P_{g,n}$).
\item [\rm(b)] If $P(S^1\times\{pt\})$ is null-homotopic, 
then $\genus(F)=1$. Moreover, $S^1\times \P(M^3)=F\times B$ and
$P\circ h=pr_2$ is the projection onto the second factor, and the fibration
$p=pr_1\co \P(M^3)=S^1\times B\rightarrow S^1$ is the projection onto the first
factor.
\end{itemize}
\end{thm}

\proof 
By lemma 3.4, the Lefschetz fibration $P$ has no
singular fibers. 

Let's first consider the case when the base $B=S^2$.
The fiber $F$ must be $T^2$, and $X$ is simply the product $F\times B$ and
$P=pr_2$ is the projection onto the second factor. This belongs to the case 
when
$P(S^1\times \{pt\})$ is null-homotopic. It is easily seen that the 
Poincar\'{e} associate $\P(M^3)$ of $M^3$ is homeomorphic to $S^1\times B$.

For the rest of the proof, we assume $\genus(B)\geq 1$.
Denote $x_0$ a base point in $X$ and let $b_0=P(x_0)$.
Consider the exact sequence induced by the fibration:
$$
1\rightarrow\pf\hookrightarrow\px\stackrel{P_*}{\rightarrow}\pb
\rightarrow 1
$$
 
Observe that $\px=\Z \oplus \pim$. Denote the generator of the $\Z$--summand
by $u=[S^1\times \{pt\}]\in\px$. 

\medskip

\noindent{\bf Case 1}\qua 
The image $P_*(u)\neq 1$ in $\pb$. 

The element $u$ is central in $\px$, therefore $P_*(u)$ is central in 
$\pb$. It follows immediately that $B$ is a torus. 
Let $n$ be the divisibility of $P_*(u)$, 
$v=\frac{1}{n}P_*(u)$
and $q$ be a simple closed loop in $B$ containing $b_0$ and representing
$v$. The monodromy $\Mon_P({P_*(u))}$ is trivial in the mapping class 
group of $F_{b_0}$. 
This is because the mapping class group of $F_{b_0}$ is isomorphic to 
$\Out(\pf)$ and $\Mon_P({P_*(u)})(\gamma)=u\gamma u^{-1}=\gamma$
for every $\gamma$ in $\pf\subset\px$ (note that $u$ is central in $\px$).

The cases when the fiber $F$ is a torus or a sphere need special treatment, and will 
be considered after the general case.

We assume for now that $\genus(F)\geq 2$.
Observe that $\Mon_P(v)$ has finite order in $\Map(F_{b_0})$.

\begin{claim}
$\Mon_P(v)=\Mon_P(q)$ has a
representative $\phi\co F_{b_0}\rightarrow F_{b_0}$ which has order $n$ and is 
periodic-point-free.
\end{claim}

\proof
Let $\phi$ be a finite order representative of $\Mon_P(q)$ provided by 
proposition 3.1.
Suppose it has periodic points. We can go back to the beginning of the proof 
and make sure
that $x_0\in F_{b_0}$ is a periodic point of $\phi$, ie,
$\phi^{n'}(x_0)=x_0$ for some $n'$ such that $0<n'<n$.
The minimal period $n'$ divides $n$ and we set $k=n/n'$.

The idea of the proof, which follows, is that if $\phi^{n'}$ has a fixed 
point, then $\phi^{n'}_*$ acts (periodically) on $\pf$ ``on a nose'',
not up to inner automorphisms of $\pf$ as in the case of periodic-point-free 
actions. This will contradict the structure of $\px$ seen from the exact 
sequence induced by the fibration.

Let $Z = F_{b_0}\times I/[_{(x,1)\sim(\phi(x),0)}]$ 
be the mapping torus of $\phi$ and 
$\widetilde{Z} = F_{b_0}\times I/[_{(x,1)\sim(\phi^{n'}(x),0)}]$ 
be the mapping torus of 
$\phi^{n'}$. 
There are a natural $n'$--fold covering $cov\co \widetilde{Z}\rightarrow Z$ and 
an embedding
$emb\co Z\hookrightarrow X$, such that $cov(x_0,0) = (x_0,0)$ and 
$ emb(x_0, 0)=x_0$.
Their composition induces a monomorphism in the fundamental groups
$(emb\circ cov)_*\co \pi_1(\widetilde{Z},(x_0,0))\hookrightarrow \px$.
Let $\delta$ be the image of 
$[\{x_0\}\times I] \in \pi_1(\widetilde{Z},(x_0,0))$ under the above monomorphism.
For any $\gamma \in \pf\subset\px$, the following relations hold in $\px$:
\begin{eqnarray*}
\delta\gamma\delta^{-1} = \phi^{n'}_*(\gamma)&\\
\delta^k\gamma\delta^{-k} = \gamma&
\end{eqnarray*}
On the other hand, we have $\delta^k=\xi u$ for some element $\xi$ in $\pf$, 
since $P_*(\delta^k)=P_*(u)$.
The relations above then imply that $\xi\gamma\xi^{-1}=\gamma$
for any $\gamma\in\pf$, thus $\xi = 1$.
Thus we have $\delta^k=u$ for $k>1$. This contradicts the fact that 
$u$ generates a direct summand in $\px$.
\endproof

Denote $\Sigma = F_{b_0}/\phi$.
Let $r$ be an oriented simple closed loop in $B$ 
such that $q\cap r = b_0$ and $\langle q,r\rangle = 1$.
The monodromy $\Mon_P(r)$ commutes with $\Mon_P(q)$. 
Hence it has a representative
$\tilde{\psi}\co F\rightarrow F$, 
which commutes with $\phi$ by lemma 3.1, 
and therefore descends to a map
$\psi\co \Sigma\rightarrow\Sigma$.
Let $\tilde{g}\co \Sigma\rightarrow S^1$ be the $S^1$--valued function
corresponding to the cyclic covering $pr\co F_{b_0}\rightarrow\Sigma= 
F_{b_0}/\phi$ and the generator $\phi$
of the deck transformations of $pr$ (cf lemma 3.3).
Since $\phi$ commutes with $\psi$, 
the function $\tilde{g}$ extends to a function 
$g\co \widetilde{M}^3\rightarrow S^1$, where $\widetilde{M}^3$ is the mapping
torus of $\psi$. In general $g$ is only a continuous function.
We can always deform it into a smooth function, which is still denoted by
$g$ for simplicity. Identify $B$ with $S^1\times S^1$
so that $(q,r)$ becomes $(S^1\times\{1\},\{1\}\times S^1)$,
then the fibrations $\pgn\co S^1\times\widetilde{M}^3\rightarrow B$
and $P\co X\rightarrow B$ have the same monodromy and therefore
$X=S^1\times M^3$ is diffeomorphic to $S^1\times \widetilde{M}^3$
by lemma 3.2.

Now we consider the case when $\genus(F)=1$.

\begin{claim}
$\Mon_P(v)=\Mon_P(q)$ is trivial.
\end{claim}

\proof
Suppose $\Mon_P(v)$ is not trivial in $\Map(F_{b_0})$.
Then we must have $n\neq 1$.
We identify $\Map(F_{b_0})$ with $SL(2,\Z)$
by the action of $\Map(F_{b_0})$ on $H_1(F_{b_0};\Z)$, and denote 
by $A\in SL(2,\Z)$ the 
element corresponding to $\Mon_P(v)$.
It has two complex eigenvalues, neither of which is equal to one.
This implies, in particular, that $\Id-A$ has non-zero determinant and 
$Q=H_1(F_{b_0};\Z)/\Im(\Id-A)$ is a finite group.
Consider the 3--manifold $Z=P^{-1}(q)$. 
It is diffeomorphic to the mapping torus of $\Mon_P(v)$.
The natural embedding of $Z$ into $X$ induces a monomorphism of fundamental groups
and $u$ is in the image of this monomorphism; thus we can regard $u$ as an 
element of $\pi_1(Z,x_0)$.
Let us recall that
$\px$ splits into a direct sum, $\px = \Z\langle u\rangle\oplus\pim$.
Since $\Z\langle u\rangle \subset \pi_1(Z,x_0)$, the fundamental group of $Z$
also splits into a direct sum, 
$\pi_1(Z,x_0) = \Z\langle u\rangle\oplus[\pi_1(Z,x_0)\cap\pim]$.
In particular, the image $[u]$ of $u$ in $H_1(Z;\Z)$ has infinite order.
Choose a loop $\delta\in \pi_1(Z,x_0)$ such that $P_*(\delta)=v$. Then
$u$ and $\delta$ are related by equation $\delta^n=u\gamma$ for
some $\gamma\in\pf$.
Denote by $[\delta]$ the image of $\delta$ in $H_1(Z;\Z)$.
Calculating homology of $Z$ we have
\[
H_1(Z;\Z)=\Z\langle[\delta]\rangle\oplus Q\text{\ \ and\ \ }
[u]=n[\delta]-[\gamma]
\]
where $Q=H_1(F_{b_0};\Z)/\Im(\Id-A)$, and 
$[\gamma]\in Q$ is the image of $\gamma$.
On the other hand, we have 
\[
H_1(Z;\Z)=\Z\langle[u]\rangle\oplus Q^\prime \text{\ \ and\ \ }
[\delta]=m[u]+[\gamma^\prime]
\]
for some integer $m\neq 0$ and $[\gamma^\prime]\in Q^\prime$,
where $Q^\prime$ is the abelianization of $\pi_1(Z,x_0)\cap\pim$.
Note that $Q^\prime$ is a finite group because 
the rank of $H_1(Z;\Z)$ is 1.
Putting the two equations together, we have
\[
(nm-1)[u]=[\gamma]-n[\gamma^\prime],\;n>1.
\]
This is a contradiction, since the left-hand-side 
has infinite order in $H_1(Z;\Z)$,
but the right-hand-side has finite order.
\endproof

Thus we have established that $\Mon_P(v)$ is trivial
and, in fact, the manifold $Z$ is diffeomorphic to $S^1\times F_{b_0}$.
It will be convenient to fix a product structure $S^1\times F_{b_0}$
on $Z$ and a flat metric on $F_{b_0}$.
These give rise to a flat metric on $Z=S^1\times F_{b_0}$.

Let us recall that $\Diff_0(F_{b_0})$ is homotopy equivalent to 
the set of conformal isometries of $F_{b_0}$; hence every 
element in $\pi_1(\Diff_0(F_{b_0}),\Id)$ could be represented by 
a linear family of parallel transforms of $F_{b_0}$.

Let $\psi\co F_{b_0}\rightarrow F_{b_0}$ be a 
representation of the monodromy of $P$ along $r$, where
$r$ is an oriented simple closed loop in $B$ 
such that $q\cap r = b_0$ and $\langle q,r\rangle = 1$.
We may assume that $\psi$ is linear with respect to the 
chosen flat metric on $F_{b_0}$.

Since $X$ fibers over a circle with fiber $Z$, 
it is diffeomorphic to the mapping torus of some self-diffeomorphism
$\Psi\co Z\rightarrow Z$, which could be chosen so that it preserves each 
fiber $\{pt\}\times F_{b_0}$.
Such a $\Psi$ is the composition of
$\Id\times \psi$ with a ``Dehn twist'' of $Z$
along $F_{b_0}$, which is defined as follows:
Choose an element in $\pi_1(\Diff_0(F_{b_0}),\Id)$ and represent it 
by a loop $\alpha\co S^1\rightarrow\Diff_0(F_{b_0})$.
Define a Dehn twist $t_{F_{b_0},\alpha}\co Z=S^1\times F_{b_0}\rightarrow Z$ by 
$t_{F_{b_0},\alpha}(s,x)=(s,\alpha(s)(x))$.
If we choose $\alpha$ to be a linear loop in the space of conformal
automorphisms of $F_{b_0}$,
then $\Psi$ will be a linear self-diffeomorphism of $Z$.

Recall that $X=S^1\times M^3$ and denote by
$\pr_1\co X\rightarrow S^1$ the projection 
onto the first factor. Let $\lambda\co Z\rightarrow S^1$ be a linear map
representing the homotopy class of $pr_1|_{Z}$,
which defines a trivial fibration on $Z$ with fiber $\Sigma\diffeo T^2$.
Then $\lambda$ and $\lambda\circ\Psi$ are homotopic,
because $\pr_1$ is defined on the whole $X$, which is the mapping torus of 
$\Psi$. This implies that $\lambda=\lambda\circ\Psi$ since both are linear
maps.

Thus $\Psi$ preserves the fibration structure of $Z$ by $\lambda$,
so that $\lambda$ extends to a fibration $\Lambda\co X\rightarrow S^1$.
Let $\phi\co \Sigma\rightarrow\Sigma$ be the self-homeomorphism induced by
$\Psi$, and denote the mapping torus of $\phi$ by $\widetilde{M}^3$
and the corresponding fibration by $p\co \widetilde{M}^3\rightarrow S^1$.
We claim that $\Lambda$ is a trivial fibration over $S^1$ with fiber
$\widetilde{M}^3$, therefore $X\diffeo S^1\times \widetilde{M}^3$. This can be
seen as follows: The monodromy $\Psi$ of $\Lambda$ is a composition of
$\Id\times \phi$ with a Dehn twist of $Z$ along $\Sigma$ since $\Psi$
preserves the trivial fibration on $Z$ by $\lambda$. The Dehn twist must
be trivial since otherwise we would have $\Psi_\ast(u)\neq u$, which 
contradicts the fact that $u$ is central in $\px$. 
Therefore $\Lambda$ is trivial
and $X=S^1\times \widetilde{M}^3$, where $\widetilde{M}^3$ is fibered
with fibration $p\co \widetilde{M}^3\rightarrow S^1$ and 
fiber $\Sigma\diffeo T^2$.
Note that the product structure $h\co X\rightarrow S^1\times \widetilde{M}^3$
could be chosen to be homotopic to the given product structure
$X=S^1\times M^3$ in the sense that $h$ is homotopic to a product of homotopy 
equivalences from $S^1$ to $S^1$ and from $M^3$ to $\widetilde{M}^3$.
Specifically, this could be done as follows: Choose a product structure
$S^1\times\Sigma$ on $Z$ by choosing a projection from $Z$ to $\Sigma$
so that the induced map on fundamental groups sends $u$ to zero.
Since both $u$ and $\Sigma$ are preserved by $\Psi$, the product structure on 
$Z$ extends to a product structure on $X$, which has the required property.

Now we will find $g\co \widetilde{M}^3\rightarrow S^1$ and $n$ such that
the fibration $P_{g,n}$ is equivalent to the original fibration $P$ on $X$.
Let $pr\co F_{b_0}\rightarrow\Sigma$ be the restriction to $F_{b_0}$ 
of the projection from $Z$ to $\Sigma$. Recall that both $F_{b_0}$ and 
$\Sigma$ are linear subspaces in $Z$. It follows that $pr$ is a 
cyclic covering
with deck transformations being parallel transforms on $F_{b_0}$.
Set $n$ to be equal to the degree of $\pr$.
Now let $\tilde{g} = (\lambda|_{F_{b_0}})^n$, where the power is 
taken point-wise in $S^1$. The function $\tilde{g}$ descends to a function 
$\tilde{\tilde{g}}\co \Sigma\rightarrow S^1$,
which is linear and $\phi$--invariant, therefore extends to a map
$g\co \widetilde{M}^3\rightarrow S^1$. It is left as an exercise for the reader
to show that $P_{g,n}$ is equivalent to $P\co X\rightarrow B$.

It remains to show that in both cases 
$\widetilde{M}^3$ is homeomorphic to the Poincar\'{e} associate $\P(M^3)$
of $M^3$. This follows from the fact that
the diffeomorphism between $S^1\times\widetilde{M}^3$ and $X=S^1\times M^3$
induces a homotopy equivalence $\widetilde{M}^3 \rightarrow M^3$,
which, by a theorem of Stallings \cite{St}, implies 
$\widetilde{M}^3 \diffeo \P(M^3)$.

When $\genus(F)=0$, $X$ is diffeomorphic to $F\times B$ 
(since $X$ is spin) with $P=pr_2$ 
the projection onto the second factor. The Poincar\'{e} associate $\P(M^3)$
is homeomorphic to $S^1\times S^2$.

\medskip

\noindent{\bf Case 2}\qua The image $P_\ast(u)=1$ in $\pb$. 

Thus $u$ lies in $\pf$ and generates a direct summand in $\pf$.
Hence the fiber $F$ must be a torus, and $u$ is primitive in $\pf$.

Identify $F$ with $S^1\times S^1$ such that the loop $q=S^1\times\{pt\}$ 
represents the class $u$ in $\pf$ and the loop $r=\{pt\}\times S^1$
represents a class in $\pi_1(M^3)$,
which we denote by $[r]$.  Then we have a reduced exact sequence:
$$
1\rightarrow \Z\langle [r]\rangle\hookrightarrow\pi_1(M^3)\stackrel
{P_\ast|_{\pi_1(M^3)}}{\longrightarrow}\pb\rightarrow 1
$$
Let $f\co \P(M^3)\rightarrow M^3$ be a homotopy equivalence between the 
Poincar\'{e} associate $\P(M^3)$ and $M^3$. It is easily
seen that there is a commutative diagram
$$
\begin{array}{ccccccccc}
1 & \rightarrow & \Z & \stackrel{j_\ast}{\rightarrow} & \pi_1(
\P({M}^3))
 & \stackrel{\pi_\ast}{\longrightarrow} & \pb & \rightarrow 
& 1\\
  &    & \| &      & \downarrow {\scriptstyle f_\ast} &   & \| &  & \\
1 & \rightarrow & \Z\langle [r]\rangle & \hookrightarrow  & \pi_1(M^3)
 & \stackrel{P_\ast|_{\pi_1(M^3)}}{\longrightarrow} & \pb & \rightarrow 
& 1.\\
\end{array}
$$

By Theorem 11.10 in \cite{H}, the Poincar\'{e} associate $\P({M}^3)$ is 
an $S^1$--fibration over the Riemann surface $B$, 
$S^1\stackrel{j}{\hookrightarrow}
\P(M^3)\stackrel{\pi}{\rightarrow}B$, from which the exact sequence
$$
1\rightarrow \Z\stackrel{j_\ast}{\rightarrow}\pi_1(\P({M}^3))\stackrel
{\pi_\ast}{\rightarrow}\pb\rightarrow 1
$$
is induced. We claim that $\P({M}^3)$
must be the trivial fibration $S^1\times B$. Suppose it is not a trivial 
fibration,
then the homology class of the fiber $j_\ast[S^1]$ is torsion in $H_1(
\P({M}^3);\Z)$, so is the fiber of the associated $T^2$--fibration of
$S^1\times \P({M}^3)$ obtained by taking the product with $S^1$, 
$S^1\times S^1\stackrel{\Id\times j}{\rightarrow}S^1\times
\P({M}^3)\stackrel{\pi}{\rightarrow}B$. But 
$((\Id\times f)\circ (\Id\times j))_\ast[S^1\times S^1]$ is homologous to the
fiber $[F]$ in $H_2(X;\Z)$, which is not torsion
since $F\hookrightarrow X\stackrel{P}{\rightarrow}B$
is a symplectic Lefschetz fibration. Hence $\P(M^3)$ is the trivial fibration
$S^1\times B$.

We then have the following commutative diagram
$$
\begin{array}{ccccccccc}
1 & \rightarrow & \Z\oplus\Z & \stackrel{(\Id\times j)_\ast}{\rightarrow} & 
\Z\oplus\Z\oplus\pb
 & \stackrel{\pi_\ast}{\longrightarrow} & \pb & \rightarrow 
& 1\\
  &    & \| &      & \downarrow {\scriptstyle(\Id\times f)_\ast} &   & \| &  & \\
1 & \rightarrow & \pf & \hookrightarrow  & \px
 & \stackrel{P_\ast}{\longrightarrow} & \pb & \rightarrow 
& 1,\\
\end{array}
$$
from which it follows that $P$ has trivial monodromy.
Let $B_0$ be a surface with boundary obtained by removing from $B$ 
a small disc $D_0$ disjoint from the base point $b_0\in\partial B_0$.
The restrictions $P|_{B_0}$ and $P|_{D_0}$ are trivial
fibrations, and $P$ is
determined by the homotopy class $[\Theta]$ of a gluing map 
$\Theta\co S^1\times F_{b_0}\rightarrow S^1\times F_{b_0}$ viewed as an
element in $Map(S^1,\Diff_0(F_{b_0}))$. 
By proposition 3.2, we have
$\pi_1(\Diff_0(F_{b_0}),\Id)=\Z\oplus\Z$. So $P$ is equivalent to a fibration
which is the product of an $S^1$--fibration over $B$, say $\xi$, with $S^1$, where
the first Chern number of $\xi$ is the divisibility of $[\Theta]$ in
$\pi_1(\Diff_0(F_{b_0}),\Id)=\Z\oplus\Z$. On the other hand, $[F]\neq 0$
in $H_2(X;\R)$, so we must have $c_1(\xi)=0$ so that 
$\Theta\sim \Id$. Thus we have proved that $X$ is
diffeomorphic to $F\times B$ with $P=pr_2$ the projection onto the second
factor.

We finish the proof of the theorem by showing the uniqueness of the symplectic structure: 
The existence of diffeomorphism
$h\co S^1\times \P(M^3)\rightarrow X$ is clear from the above classification
of the Lefschetz fibration $P$. Moreover, the canonical symplectic form
on $S^1\times \P(M^3)$ is positive on each fiber of the pull-back fibration
$P\circ h$, hence by lemma 2.1, it is deformation equivalent to
$h^\ast\omega$.
\endproof

\end{document}